\documentclass[11pt]{article}
\usepackage{geometry}                
\usepackage{graphicx}
\usepackage{amssymb}
\usepackage{epstopdf}
\usepackage{amsmath,amsthm,amscd,amssymb}
\numberwithin{equation}{section}

\theoremstyle{plain}
\newtheorem{theorem}{Theorem}[section]
\newtheorem{lemma}[theorem]{Lemma}
\newtheorem{corollary}[theorem]{Corollary}

\newtheorem{conjecture}[theorem]{Conjecture}

\theoremstyle{definition}
\newtheorem{definition}[theorem]{Definition}

\theoremstyle{remark}
\newtheorem{remark}[theorem]{Remark}

\newtheorem{case[theorem]}{Case}

\title{Averages over hyperplanes, sum-product theory in vector spaces
over finite fields and the Erd\H os-Falconer distance conjecture}
\author{Derrick Hart, Alex Iosevich, Doowon Koh and Misha Rudnev}

\begin{document}

\maketitle

\begin{abstract} We prove a point-wise and average bound for the
number of incidences between points and hyper-planes in vector spaces
over finite fields. While our estimates are, in general, sharp, we
observe an improvement for product sets and sets contained in a
sphere. We use these incidence bounds to obtain significant
improvements on the arithmetic problem of covering ${\mathbb F}_q$,
the finite field with $q$ elements, by $A \cdot A+\dots +A \cdot A$,
where $A$ is a subset ${\mathbb F}_q$ of sufficiently large size. We
also use the incidence machinery we develope and arithmetic
constructions to study the Erd\H os-Falconer distance conjecture in
vector spaces over finite fields.  We prove that the natural analog of
the Euclidean Erd\H os-Falconer distance conjecture does not hold in
this setting due to the influence of the arithmetic.  On the positive
side, we obtain good exponents for the Erd\H os -Falconer distance
problem for subsets of the unit sphere in $\mathbb F_q^d$ and discuss
their sharpness. This results in a reasonably complete description of
the Erd\H os-Falconer distance problem in higher dimensional vector
spaces over general finite fields.
 \end{abstract}

\tableofcontents

\section{Introduction}
Let  ${\mathbb F}_q$ denote a finite field with $q$ elements, where
$q$, a power of an odd prime, is viewed as an asymptotic parameter. In
a special case when $q=p$ is a prime,  we use the notation
${\mathbb Z}_p$. Let ${\mathbb F}_q^*$ denote the multiplicative group
of ${\mathbb F}_q$. How large does $A \subset {\Bbb F}_q$ need to be
to make sure that
$$ dA^2=\underbrace{A^2+\dots+A^2}_{d \mbox{ {\small times}
}}\supseteq{\mathbb F}^*_q?$$

Define
$$ A^2=A\cdot A= \{a \cdot a': a,a' \in A\}\;\mbox{ and }\;A+A=\{a +
a': a,a' \in A\}.$$

It is known (see e.g. \cite{GK06}) that if $d=3$ and $q$ is prime,
this conclusion holds if the number of elements $|A| \ge
Cq^{\frac{3}{4}}$, with a sufficiently large constant $C>0$. It is
reasonable to conjecture that if $|A| \ge
C_{\epsilon}q^{\frac{1}{2}+\epsilon}$, then $2A^2\supseteq{\mathbb
F}_q^*$. This
result cannot hold, especially in the setting of general finite fields
if $|A|=\sqrt{q}$ because $A$ may in fact be a subfield. See also
\cite{BGK06}, \cite{C04}, \cite{G06}, \cite{G07}, \cite{HIS07},
\cite{KS07}, \cite{TV06}, \cite{V07} and the references contained
therein on recent progress related to this problem and its analogs.
For example, Glibichuk, \cite{G06}, proved that
$$ 8A \cdot B={\mathbb Z}_p, $$ $p$ prime, provided that $|A||B|>p$ and
either $A=-A$ or $A \cap (-A)=\emptyset$. Glibichuk and Konyagin,
\cite{GK06}, proved that if $A$ is subgroup of ${\mathbb Z}_p^{*}$,
and
$|A|>p^{\delta}$, $\delta>0$, then
$$ NA={\mathbb Z}_p$$ with
$$ N \ge C4^{\frac{1}{\delta}}.$$ The above-mentioned results were
achieved by methods of arithmetic combinatorics.

One of the goals of this paper is to use the geometry of the vector
space ${\mathbb F}_q^d$, where $q$ is not necessarily a prime number,
to deduce a good lower bound on the size of $A$ that guarantees that
$dA^2\supseteq{\mathbb F}^*_q$.

The second aim of this paper is directly related to the finite field
version of the Erd\H os-Falconer distance problem. The Erd\H os
distance conjecture says that if $E$ is a finite subset of ${\mathbb
R}^d$, $d \ge 2$, then
$|\Delta(E)| \ge C_{\epsilon}{|E|}^{\frac{2}{d}-\epsilon},$ where
$\Delta(E)=\{\|x-y\|:x,y \in E\}$,
$\|\cdot\|$ denotes the standard Euclidean metric. This problem is far
from resolution in any
dimension. See,
for example, a monograph by Matousek (\cite{Ma02}) and the
references contained therein to review the main milestones of the
progress towards this conjecture.

The Falconer distance conjecture says that if $E \subset {\mathbb R}^d$,
$d \ge 2$, has Hausdorff dimension greater than $\frac{d}{2}$, then
$\Delta(E)$ has positive Lebesgue measure. See
\cite{Erd05} for the latest progress and description of techniques.
For the connections between the Erd\H os and Falconer distance
problems see, for example, \cite{IRU07}.

In the finite field setting the question turns out to have features of
both the Erd\H os and Falconer distance problems. The first
non-trivial result was obtained by Bourgain, Katz and Tao (\cite{BKT04})
using
arithmetic-combinatorial methods and the connection of the geometric
incidence problem of counting distances with sum-product estimates.
\begin{theorem} \label{bkt} Suppose $E \subset {\mathbb Z}_p^2$, where $p
\equiv 3 \mod 4$ is a prime, and $|E| \leq p^{2-\epsilon}$. Then
there exists $\delta=\delta(\epsilon)$ such that
$$|\Delta(E)| \geq c{|E|}^{\frac{1}{2}+\delta}.$$
\end{theorem}
Here and  throughout the paper, for $E \subseteq {\mathbb F}_q^d$,
$$\Delta(E)=\{\|x-y\|={(x_1-y_1)}^2+\dots+{(x_d-y_d)}^2: x,y \in E\}$$
denotes the distance set of $E$.

It is interesting to observe that while the quantity $\|\cdot \|$ is
not a distance, in the traditional sense, it is still a natural object
in that it is invariant under the action of orthogonal matrices.

We note that the conclusion of Theorem \ref{bkt} with the exponent
$\frac{1}{2}$ follows from
the argument due to Erd\H os (\cite{Erd45}). The condition $|E|
\lesssim q^{2-\epsilon}$ addresses
the fact that if $E={\mathbb Z}_p^2$, then $\Delta(E)={\mathbb Z}_p$
and so $|\Delta(E)|=\sqrt{|E|}$ and no better. The condition $p \equiv
3 \mod 4$ addresses the fact that if conversely $p \equiv 1 \mod 4$,
the field ${\mathbb F}_p$ contains an element $i$ such that $i^2=-1$.
This would allow one to take
\begin{equation}\label{nl}E=\{(t,it): t \in {\mathbb
Z}_p\}\end{equation} and it is straightforward to check that while
$|E|=p$, $|\Delta(E)|=1$ as all the distances between the elements of
the set are identically
$0$.

In view of the examples cited in the previous paragraph, Iosevich and
Rudnev (\cite{IR07}) formulated the Erd\H os-Falconer conjecture as
follows.
\begin{conjecture} \label{kingofsweden} Let $E \subset {\mathbb F}_q^d$
such that $|E| \ge C_{\epsilon}q^{\frac{d}{2}+\epsilon}$. Then there
exists $c>0$ such that
$$ |\Delta(E)| \ge cq.$$
\end{conjecture}

A Fourier analytic approach to this problem, developed in \cite{IR07},
led to the following result.
\begin{theorem} \label{ir} Suppose that $E \subset {\mathbb F}_q^d$ and
$|E| \ge 4q^{\frac{d+1}{2}}$. Then $\Delta(E)={\mathbb F}_q$.
\end{theorem}

In $d=2$, in particular, whenever $q\ll|E|\ll q^2$, the results of the
paper provide a quantitatively
explicit estimate. This enabled the first two authors and J. Solymosi
to use its analog to obtain a strongly nontrivial sum-product estimate
(\cite{HIS07}). The basis for the latter results was Weil's
bound (\cite{We48}) for Kloosterman sums,
$$ \left| \sum_{t \not=0} \chi(at+t^{-1}) \right| \leq 2\sqrt{q},$$
where $\chi$ further denotes a non-trivial additive character of
${\mathbb F}_q$.

Observe that in the formulation of Conjecture \ref{kingofsweden} one
asks for the {\em positive proportion} of distances in ${\mathbb
F}_q$, while Theorem \ref{ir} guarantees that {\em all} distances in
${\mathbb F}_q$ occur, being generated by $E$. The latter question is
closely related to what in the discrete Euclidean setting is known as
the Erd\H os {\em single distance} conjecture, which says that a
single distance in $\mathbb R^2$ cannot occur more than
 $c n^{1+\epsilon}$ times where $n$ is the cardinality of the
underlying point set $E$. It is tempting to strengthen the claim of
Conjecture \ref{kingofsweden} to cover all distances. However, we
shall see below that even the weak form of this conjecture
(\ref{kingofsweden}) is not true.
This shows that the Theorem (\ref{ir}) is essentially sharp.  This
underlines the difference between the finite field setting and the
Euclidean setting where the Erd\H os-Falconer distance conjecture,
while far from being proved, is still strongly believed.

We shall see however, that the exponent predicted by Conjecture
\ref{kingofsweden} does hold for subsets of the sphere in
$$S^{d-1}=S=\{x \in {\mathbb F}_q^d: x_1^2+\dots+x_d^2=1\}$$
in even dimensions. While it is possible that in some cases this
exponent may be further improved for this class of sets under some
circumstances, we provide examples showing that if one is after all
the distances, and not a positive proportion, then getting a better
index is not in general possible. This is geometrically analogous to
the general case, for since $S^{d-1}$ is $(d-1)$-dimensional variety
in
${\Bbb F}_q^d$, it makes sense that the sharp index should be
$\frac{(d-1)+1}{2}=\frac{d}{2}$. The motivation for studying the Erd\H
o-Falconer distance problems for subsets of the sphere is not limited
by the consideration that it provides a large set of sets for which
Conjecture \ref{kingofsweden} holds. For example, Erd\H os original
argument that shows that $N$ points in ${\Bbb R}^2$ determine $\gtrsim
N^{\frac{1}{2}}$ distances proceeded as follows. Choose one of the
points in the set and draw circles of every possible radius centered
at this point such that each circle contains at least one other point
of the set. Suppose that the number of such circles is $t$. If $t \ge
N^{\frac{1}{2}}$, we are done. If not, there exists a circle
containing $\ge N/t$ points and these points, by an elementary
argument, determine $\ge N/2t$ distinct distances. Comparing $t$ and
$N/2t$ yields the conclusion. In higher dimension we may proceed by
induction with the induction hypothesis being the number of distances
determined by points on a sphere. Thus one may view the distribution
of distances determined by points on a sphere as a natural and
integral component of the general Erd\H os distance problem.

We conclude our introduction by emphasizing that the proofs below show
that the arithmetic structure of general fields allows for example
that may not have analogs in Euclidean space. A detailed comparative
study between the Euclidean and finite field environments shall be
conducted in a subsequent paper.

\vskip.125in

\section{Statement of results}

\subsection{Key incidence estimate} Our main tool is the following
incidence theorem. See \cite{HI07} for an earlier version.
\begin{theorem} \label{geom} Let $E \subset {\mathbb F}_q^d$ and
define the {\rm incidence function}
\begin{equation}\label{if} \nu(t)=\{(x,y) \in E \times E: x \cdot
y=t\}.\end{equation}

Then
\begin{equation} \label{L2} \sum_{t \in {\mathbb F}_q }\nu^2(t) \leq
{|E|}^4q^{-1}+|E|q^{2d-1} \sum_{k \neq (0,\dots,0)}
|E \cap l_k|{|\widehat{E}(k)|}^2+(q-1)q^{-1}{|E|}^2E(0, \dots, 0),
\end{equation} where
\begin{equation}\label{linek} l_k=\{tk: t \in {\mathbb
F}^*_q\}.\end{equation}

Moreover,
\begin{equation} \label{pointwise} \nu(t)={|E|}^2q^{-1}+R(t),
\end{equation} with
\begin{equation}\hspace{.15in}\left\{\begin{array}{llllll} |R(t)|
&\leq &|E|q^{\frac{d-1}{2}}, &\mbox{ for }\;t \not=0, \\ \hfill \\
|R(0)| &\leq &|E|q^{\frac{d}{2}}.\end{array}\right.\label{trb}\end{equation}
\end{theorem}

Note that $E(x)$ denotes the characteristic function of $E$, so $E(0,
\dots, 0)=1$ if the origin is in $E$ and $0$ otherwise. Also note that
in many of the applications below it is legitimate to assume, without
loss of generality, that $E$ does not in fact contain the origin.

\begin{remark} The proof of Theorem \ref{geom} is via Fourier
analysis. It has been pointed out to the authors by Seva Lev that an
alternate approach to (\ref{pointwise}) is via a graph theoretic
result due to Alon and Krivelevich. See, \cite{AK97} and the
references contained therein. We also note that the relevant result of
Alon and Krilevich can be recovered from the estimate (\ref{trb})
above.
\end{remark}

\begin{remark} \label{mattila} There are parallels here that are worth
pointing out. In the study of the Euclidean Falconer conjecture, the
$L^2$ norm of the distance measure is dominated by the Mattila
integral, discovered by P. Mattila, (\cite{Mat87}):
$$ \int_1^{\infty} {\left( \int_{S^{d-1}} {|\widehat{\mu}(t\omega)|}^2
d\omega \right)}^2 t^{d-1}dt,$$ where
$\mu$ is a Borel measure on the set $E$ whose distance set is being
examined. It is reasonable to view the expression
$$ \sum_{k \neq (0,\dots,0)} |E \cap l_k| {|\widehat{E}(k)|}^2$$ as
the Mattila integral for the dot product problem, a direct analog of
the Mattila integral for the distance set problem in the Euclidean
space.
\end{remark}
By analogy with the distance set $\Delta(E)$, let us introduce the set
of dot products
\begin{equation}
\Pi(E)=\{x\cdot y=x_1y_1+\ldots+x_d y_d:\,x,y\in E\},
\label{Pi}\end{equation}

\begin{corollary} \label{L2dot} Let $E \subset {\mathbb F}_q^d$ such that
$|E|>q^{\frac{d+1}{2}}$. Then
$$ {\mathbb F}_q^{*} \subseteq \Pi(E).$$

This result cannot in general be improved in the following sense:
\renewcommand{\theenumi}{\roman{enumi}}
\begin{enumerate}  \item Whenever ${\mathbb F}_q$ is a quadratic
extension, for any
$\epsilon>0$ there exists $E \subset {\mathbb F}_q^d$ of size $
\approx q^{\frac{d+1}{2}-\epsilon}$, such that $|\Pi(E)|=o(q)$. In
particular, the set of dot products does not contain a positive
proportion of the elements of ${\mathbb F}_q$. \item For
$d=4m+3,\,m\geq0$, for any $q\gg1$ and any $t\in {\mathbb F}_q^*$,
there exists $E$ of cardinality $ \approx q^{\frac{d+1}{2}}$, such
that $t\not\in\Pi(E).$ \end{enumerate}
\end{corollary}

Throughout the paper, $X \lesssim Y$ means that there exists $C>0$
such that $X \leq CY$.

\vskip.125in

\subsection{Arithmetic results}
\begin{theorem} \label{kickass} Let $A \subset {\mathbb F}_q$, where
${\mathbb
F}_q$ is an arbitrary finite field with $q$ elements, such that
$|A|>q^{\frac{1}{2}+\frac{1}{2d}}$. Then
\begin{equation} \label{sex} {\mathbb F}_q^{*} \subset dA^2.\end{equation}

Moreover, suppose that for some constant $C^{\frac{1}{d}}_{size}$,
$$|A| \ge C^{\frac{1}{d}}_{size} q^{\frac{1}{2}+\frac{1}{2(2d-1)}}.$$
Then
\begin{equation} \label{sex2} \hspace{-5mm}|dA^2| \ge q \cdot
\frac{C^{2-\frac{1}{d}}_{size}}{C^{2-\frac{1}{d}}_{size}+1}. \end{equation}
\end{theorem}

It follows immediately from Theorem \ref{kickass} that in the most
interesting particular case $d=2$,
$$\hspace{9mm}{\mathbb F}_q^{*} \subset A^2+A^2\quad \mbox{ if }
\quad |A|>q^{\frac{3}{4}},$$ and\hspace{-8mm}
$$ |A^2+A^2| \ge q \cdot
\frac{C^{\frac{3}{2}}_{size}}{C^{\frac{3}{2}}_{size}+1}\quad \mbox{ if }\quad
|A| \ge C_{size}^{\frac{1}{2}} q^{\frac{2}{3}}.$$

\vskip.125in
We would like to complement the general result in Theorem
\ref{kickass} with the following conditional statement.
\begin{theorem} \label{uniform} Let $A \subset {\mathbb F}_q$, with
$|A| \ge C^{\frac{1}{2}}_{size}q^{\frac{1}{2}}$, and suppose that
\begin{equation} \label{otyebis} |(A \times A) \cap t(A \times A)|
\leq C_{uni} {|A|}^2q^{-1},
\end{equation} for all $t \in {\mathbb F}^*_q\setminus\{1\}$.
Then
$$ |2A^2| \ge q \cdot \frac{C_{size}}{2C_{size}+C_{uni}}.$$
\end{theorem}

\vskip.125in

\subsection{Distance set results}
\begin{theorem}\label{erdosfalse} The Conjecture \ref{kingofsweden} is
false. More precisely,
there exists $c>0$ and $E \subset \mathbb F_q^d$, $d$ is odd, such that
$$|E| \ge cq^{\frac{d+1}{2}} \quad \rm{and} \quad \Delta(E) \not={\Bbb
F}_q.$$
\end{theorem}

\renewcommand{\theenumi}{\roman{enumi}}

\begin{theorem} \label{sphericalerdos} Let $E \subset {\mathbb
F}_q^d$, $d \ge 3$, be a subset of the sphere $S=\{x \in {\mathbb
F}_q^d:\, \|x\|=1\}$.

\begin{enumerate} \item Suppose that $|E| \ge Cq^{\frac{d}{2}}$ with a
sufficiently large constant $C$. Then there exists $c>0$ such that
\begin{equation} \label{pospropsphere} |\Delta(E)| \ge cq.\end{equation}

\item If $d$ is even, then under the same assumptions as above,
\begin{equation} \label{spherepointwise} \Delta(E)={\mathbb F}_q.
\end{equation}

\item If $d$ is even, there exists $c>0$ and $E \subset S$ such that
\begin{equation} \label{evencounter}|E| \ge cq^{\frac{d}{2}} \quad
\rm{and} \quad
 \Delta(E) \not={\mathbb F}_q. \end{equation}

\item If $d$ is odd and $|E| \ge Cq^{\frac{d+1}{2}}$ with a
sufficiently large constant $C>0$, then
\begin{equation} \label{dotall} \Delta(E)={\mathbb F}_q.\end{equation}

\item If $d$ is odd, there exists $c>0$ and $E \subset S$ such that
\begin{equation} \label{oddcounter} |E| \ge cq^{\frac{d+1}{2}} \quad
\rm{and} \quad
 \Delta(E) \not={\mathbb F}_q. \end{equation}
\end{enumerate}

\end{theorem}

\begin{remark} In summary, we always get a positive proportion of all
the distances if $|E| \ge Cq^{\frac{d}{2}}$. If $d$ is even, we get
all the distances under the same assumption and the
size condition on $E$ cannot be relaxed. Similarly, if $d$ is odd we know
that we cannot in general get all the distances if
$|E|\ll q^{\frac{d+1}{2}}$, but, as we note above, we get a positive
proportion of the distances under the assumption that $|E| \ge
Cq^{\frac{d}{2}}$, and it is not out of the question that one can go
as low as $q^{\frac{d-1}{2}+\epsilon}$, asking for the positive
proportion of distances. \end{remark}

We conclude this section by formulating a result which says that if a
subset of the sphere is statistically evenly distributed, then the
distance set is large under much milder assumptions than above.
\begin{definition} \label{welldistributed} Let
$E \subset S=\{x \in {\mathbb F}_q^d: x_1^2+\dots+x_d^2=1\}$. Suppose that
$$ |E \cap H| \leq C|E|q^{-1}$$ for every $(d-1)$-dimensional
hyper-plane $H$ passing through the origin. Then we say that $E$ is
uniformly distributed on the sphere.
\end{definition}

We note that since the density of $E$ is $\frac{|E|}{q^d}$ and the
density of a hyperplane $H$ is
$\frac{|H|}{q^d}=q^{-1}$, the expected number of points on $E \cap H$
is indeed $q^d \cdot \frac{|E|}{q^d} \cdot q^{-1}=\frac{|E|}{q}$. Thus
the uniformity assumption says that the number of points of $E$ on
each hyperplane through the origin does not exceed the expected number
by more than a constant.

\begin{theorem} \label{uniformsphere} Suppose that $E$ is uniformly
distributed on the sphere and that $|E| \ge Cq$. Then
\begin{equation} \label{lafa} |\Delta(E)| \ge cq. \end{equation}
\end{theorem}

\vskip.125in

\subsection{Acknowledgements:} The authors wish to thank Luca
Brandolini, Leoardo Colzani, Giaccomo Gigante, Nets Katz, Sergei
Konyagin, Seva Lev, Igor Shparlinsky and Ignacio Uriarte-Tuero for
many helpful remarks about the content of this paper.

\vskip.125in

\section{Proof of geometric results: Theorem \ref{geom} and Corollary
\ref{L2dot}}

\vskip.125in

\subsection{Proof of  the $L^2$ estimate (\ref{L2}):}

The Fourier transform of a complex-valued function $f$ on ${\mathbb
F}^d_q$ with respect to a non-trivial principal additive character
$\chi$ on ${\Bbb F}_q$ is given by
$$\widehat f(k) = q^{-d}\sum_{x \in {\mathbb F}^d_q} \chi(-x \cdot k)
f(x)$$ and the Fourier inversion formula takes the form
$$f(x) = \sum_{k\in {\mathbb F}^d_q} \chi(x \cdot k) \widehat{f}(k)$$

\vskip.125in

We have
$$ \begin{array}{lll}\nu(t)&=&|\{(x,y) \in E^2: x \cdot y=t\}|\\ \hfill \\
&=& \sum_{x \cdot y=t} E(x)E(y).\end{array}$$

The Cauchy-Schwartz inequality applied to the sum in the variable $x$ yields
\begin{equation}\begin{array}{lll}\sum_t \nu^2(t) &\leq& |E| \cdot
\sum_t \sum_{x \cdot y=t} \sum_{x \cdot y'=t} E(x)E(y)E(y')\\ \hfill
\\ &=&|E| \sum_{(y'-y) \cdot x=0} E(y')E(y)E(x)\\ \hfill \\ &=&|E|
q^{-1} \sum_{y',y,x} \sum_s \chi(s((y'-y) \cdot x)) E(y')E(y)E(x)\\
\hfill \\&=
&{|E|}^4q^{-1}+|E|q^{2d-1} \sum_x \sum_{s \not=0} E(x)
{|\widehat{E}(sx)|}^2\\ \hfill \\ &=
&{|E|}^4q^{-1}+|E|q^{2d-1} \sum_x \sum_{s \not=0} E(sx)
{|\widehat{E}(x)|}^2\\ \hfill \\ &=
&{|E|}^4q^{-1}+|E|q^{2d-1} \sum_{x\neq(0,\ldots,0)} |E \cap l_x|
{|\widehat{E}(x)|}^2+(q-1)q^{-1}{|E|}^3E(0, \dots,
0).\end{array}\label{vl}\end{equation} In the third line we have used
the standard trick that $\sum_{s\in \mathbb F_q} \chi(ts)$ equals $q$
for $t=0$ and zero otherwise. The transition from the fourth line to
the fifth one was after changing variables $sx \to x$ and then $s \to
s^{-1}.$ This completes the proof of the estimate (\ref{L2}), which
uses the ``Fourier'' notation $k$ for $x$.

\vskip.125in

\subsection{Proof of the point-wise estimate (\ref{pointwise})}
Similarly to the third line of (\ref{vl}), we rewrite the expression
for the incidence function (\ref{if}) in the form
$$ \nu(t)=\sum_{x,y \in E} q^{-1} \sum_{s \in {\mathbb F}_q} \chi(s(x \cdot
y-t)).$$

Isolating the term $s=0$ we have, according to (\ref{pointwise})
\begin{equation}\label{break}\begin{array}{llll}
\nu(t)&=&{|E|}^2q^{-1}+R(t), \quad \mbox{where}
\\ \hfill \\ R(t)&=&\sum_{x,y \in E} q^{-1} \sum_{s \not=0} \chi(s(x
\cdot y-t)).\end{array}\end{equation}

Viewing $R$ as a sum in $x$, applying the Cauchy-Schwartz inequality and
dominating the sum over $x \in E$ by the sum over $x \in {\mathbb F}_q^d$, we
see that
$$\begin{array}{llll}R^2(t) &\leq &|E| \sum_{x \in {\mathbb F}_q^d}
q^{-2} \sum_{s,s' \not=0}
\sum_{y,y' \in E}
\chi(sx \cdot y-s'x \cdot y') \chi(t(s'-s))\\ \hfill \\ &=&|E| q^{d-2}
\sum_{\substack{sy=s'y' \\ s,s' \not=0}}
\chi(t(s'-s))E(y)E(y')\\ \hfill \\ &=&I + II,\end{array}$$ whether the
term $I$ corresponds to the case $y=y'$ (which forces $s=s'$), and the
term $II$ corresponds to the case $y\neq y'$ (and so $s\neq s'$).

In the latter case we may set $a=s/s', b=s'$ and obtain, for $t\neq0$,
\begin{equation}\begin{array}{llll} II &=& |E| q^{d-2}
\sum_{y,b\neq0;\,a \neq 0,1}
\chi(tb(1-a))E(y)E(ay)\\ \hfill \\ &=& -|E| q^{d-2} \sum_{y,a \not=1,0}
E(y)E(ay).\end{array} \label{tak}\end{equation}
Thus,
\begin{equation}
\begin{array}{llll} |II(t)|& \leq & |E|q^{d-2} \sum_{y \in
E\setminus\{(0,\ldots,0)\}} (|E \cap l_y|+1)\\ \hfill \\ &\leq&
{|E|}^2 q^{d-1},\end{array}\label{takk}\end{equation}  since
$|E \cap l_y|+1 \leq q$ by the virtue of the fact that each straight
line contains
exactly $q$ points. The term $+1$ above has been added because in
(\ref{linek}) above, the line $l_y$ was defined away from the origin.

In the case $s=s'$ we get
\begin{equation}\label{esti}I(t)= |E|q^{d-2} \sum_{s\neq 0; y} E(y) <
{|E|}^2q^{d-1}.\end{equation}
It follows that for $t\neq0$
$$ R^2(t) \leq -Q(t)+{|E|}^2q^{d-1},$$ with
$$ Q(t) \ge 0.$$
Therefore, for $t\neq0$ we have the bound (\ref{trb}),
\begin{equation} \label{remainder} |R(t)| \leq
|E|q^{\frac{d-1}{2}}.\end{equation}

The same argument shows that
$$ |R(0)| \leq |E|q^{\frac{d}{2}}.$$

\vskip.125in

\subsection{Proof of Corollary \ref{L2dot}} We now turn our attention
to the Corollary \ref{L2dot}. The sufficient condition for
$\Pi(E)\supseteq\mathbb F_q^*$ follows immediately from (\ref{break})
and (\ref{remainder}). Quite simply, it follows that $\nu(t)>0$ for
all $t\neq 0$.

\medskip
To address the statement (i) of the Corollary, let us consider the
case $d=2$ and $q=p^2$, where $p$ is a power of a large prime. The
higher dimensional case follows similarly. Let $a$ be a generator  of
the cyclic group ${\mathbb F}_q^{*}$. Then $a^{q-1}=1$ and $a^{p+1}$
is the generating element for ${\mathbb F}_p^{*}$ since
$p+1=\frac{q-1}{p-1}$.

Let $A$ be a proper cyclic subgroup of ${\mathbb F}_q^{*}$ which
properly contains ${\mathbb F}_p^{*}$. Let $s$ be a divisor of $p+1$ and
let the generating element of $A$ be $\alpha=a^s$. Note that we are
taking advantage of the fact that ${\mathbb F}_q^{*}$ is cyclic.
Consider the unit circle
$$\{x \in {\mathbb F}_q^2: x_1^2+x_2^2=1\},$$ and its subset
$$ C_p=\{x \in {\mathbb F}_p^2: x_1^2+x_2^2=1\}.$$

By elementary number theory (or Lemma \ref{sphereft}), the cardinality
of $C_p$ is $p\mp 1$, depending on whether
negative one is or is not a square in ${\mathbb F}_p^{*}$. Clearly,
for any $u,v\in C_p$,
$u \cdot v \in {\mathbb F}_p$. Let
\begin{equation} E=\{tu: t\in A, u\in C_p\}. \label{cset} \end{equation}

For any $x,y \in E$, the dot product $x\cdot y$, if nonzero, will lie
in $A$. Indeed, if $x=tu$, $y=\tau v$, according to (\ref{cset}), then
$$x\cdot y=t\tau(u\cdot v)\in A\cup\{0\},$$ since $A$ contains
${\mathbb F}_p^{*}$. The cardinality of $E$ is
$$|E| = \frac{p\mp 1}{2} |A| = \frac{p\mp 1}{2} \cdot \frac{q-1}{s},$$
where $s$ is a divisor of $p+1$. Taking
$s=2$ works and shows that less than half the elements of ${\mathbb
F}_q^{*}$ may be realized as dot products determined by a set of size
in excess of $ \frac{1}{4} \cdot {q}^{\frac{3}{2}}$. In order to see
that
$\{x \cdot y: x,y \in E\}$ does not in general even contain a positive
proportion of the elements of
${\mathbb F}_q$ if $|E|\ll q^{\frac{3}{2}}$, we need to produce a
sequence of primes, or prime powers, such that $p+1$ has large
divisors. For the reader who does not believe in the existence of
infinitely many generalized Fermat primes (those in the form
$a^{2^n}+1$), we can always do it using field extensions as follows.

Consider the family of prime powers
$$ \{p^{2k+1}: k=1,2 \dots \}$$ and observe that
$$ p+1 \ {|} \ p^{2k+1}+1.$$

This completes the construction demonstrating the claim (i). To take
care of the higher dimensional case, simply replace circles by spheres
and the same argument goes through.

\medskip
The claim (ii) of the Corollary will follow immediately from the
construction used in the proof of the item (v) of Theorem
\ref{sphericalerdos} (see Section \ref{nodot}).  On any sphere
$\{x \in {\Bbb F}_q^d: x_1^2+\dots+x_d^2=r\}$, with $d=4m+3$, we can
find a set $E$, with $|E| \gtrsim q^{\frac{d+1}{2}},$ such that the
dot product $t=-r$ is not achieved.

\vskip.125in

\section{Proof of the arithmetic results}

\subsection{Proof of Theorem \ref{kickass}}
We may assume, without loss of generality, that $A$ does not contain
$0$. Let ${\displaystyle E=\underbrace{A^2+\dots+A^2}_{d \mbox{
{\small times} }}}$. The proof of the first part of
Theorem \ref{kickass} follows instantly from the estimate
(\ref{pointwise}). To prove the second part observe that
$$ |E \cap l_y| \leq |A|={|E|}^{\frac{1}{d}}$$ for every $y \in E$.
Using this, the estimate (\ref{L2}) implies, by  the Cauchy-Schwartz
inequality
\begin{equation} \label{csagain} {|E|}^4={\left( \sum_t \nu(t)
\right)}^2 \leq |\Pi(E)| \cdot \sum_t \nu^2(t),\end{equation}
that
$$ |\{(x \cdot y: x,y \in E\}| \ge q \cdot \frac{{|E|}^2}{q^d \cdot
{|E|}^{\frac{1}{d}}+{|E|}^2},$$ and, consequently, that
$$|\{x \cdot y: x,y \in E\}| \ge q \cdot
\frac{C^{2-\frac{1}{d}}_{size}}{C^{2-\frac{1}{d}}_{size}+1} $$
if
$$|E| \ge C_{size}q^{\frac{d}{2}+\frac{d}{2(2d-1)}}.$$

It follows that if
$$ |A| \ge C^{\frac{1}{d}}_{size}q^{\frac{1}{2}+\frac{1}{2(2d-1)}},$$ then
$$ |dA^2| \ge q \cdot
\frac{C^{2-\frac{1}{d}}_{size}}{C^{2-\frac{1}{d}}_{size}+1}$$ as
desired. This
completes the proof of Theorem \ref{kickass}.

\vskip.125in

\subsection{Proof of the conditionally optimal arithmetic result (Theorem
\ref{uniform})}

\vskip.125in

Once again throw zero out of $A$ if it is there, and let $E=A \times
A$. Using (\ref{L2}) we see that
$$ \sum_t \nu^2(t) \leq {|E|}^4q^{-1}+q^3|E| \sum_{k\neq (0,0)} |E \cap l_k|
{|\widehat{E}(k)|}^2.$$
Now,
$$ \begin{array}{lll} q^3|E| \sum_{k\neq (0,0)} |E \cap l_k|
{|\widehat{E}(k)|}^2& \leq &q^{3}|E| \cdot |E|
\cdot {|\widehat{E}(1,1)|}^2\\ \hfill \\  &+& q^3|E| \sum_{k
\not=(0,0),(1,1)} |E \cap l_k| {|\widehat{E}(k)|}^2\\ \hfill \\ &\leq
&{|E|}^4q^{-1}+C_{uni}{|E|}^3.\end{array}$$

It follows that
$$ \begin{array}{lll} |2A^2|&= &|\{x \cdot y: x,y \in E\}| \\ \hfill
\\ & \ge & \frac{{|E|}^4}{{|E|}^4q^{-1}+C_{uni}{|E|}^3}\\ \hfill \\ &
\ge & q\cdot
\frac{C_{size}}{2C_{size}+C_{uni}},\end{array}$$
as desired.

\vskip.125in

\section{Distances: Proofs of Theorems \ref{erdosfalse},
\ref{sphericalerdos} and \ref{uniformsphere}}
\vskip.125in
\subsection{Proof of Theorem \ref{erdosfalse}}

We begin by proving the following lemma.

\begin{lemma}  \label{nullvectors} We say that $v \in \mathbb F_q^d$,
$v \not=(0, \dots, 0)$,
is a {\it null} vector if  $v\cdot v=0$. If $d\geq4$ is even, then
there exists $\frac{d}{2}$ mutually orthogonal null vectors
$v_1,\dots, v_{\frac{d}{2}}$ in ${\Bbb F}_q^d$.  \end{lemma}

To prove the lemma, suppose there exists an element $i \in \mathbb
F_q$ such that $i^2=-1$. Consider the collection of vectors
$$v_1=(1,i,0,0,\ldots,0,0), \ v_2=(0,0,1,i,\ldots,0,0), \dots ,
v_{\frac{d}{2}}=(0,0,\ldots,0,0,1,i).$$

It follows immediately that
$$v_k\cdot v_l=0$$ for every $k,l=1,\ldots, \frac{d}{2}$.

If $-1$ is not a square, then from simple counting there exists a null
vector
$$v_1=(a,b,c,0\ldots,0),$$ with all $a,b,c \in \mathbb F_q^*.$ Suppose, $d$ is a multiple of $4$.
Then we may take the null vector
$$v_2=(0,-c,b,a,\ldots,0)$$ noting that this vector is orthogonal to

$v_1$.  In this same way we may now take the null vector

$$v_3=(0,0,0,0,a,b,c,0,\ldots,0),$$ and find a corresponding null vector $v_4$ which is orthogonal to $v_3$ as well as trivially orthogonal to $v_1$ and $v_2$.  Continuing in this manner we obtain
$\frac{d}{2}$ mutually orthogonal null vectors.

The proof will be complete now if we can also treat the case $d=6$. In this case
Let
$$ v_1=(a,b,c,0,0,0),\;\;v_2=(0,0,0,a,b,c) \ \text{where} \ a^2+b^2+c^2=0.$$

Consider two three-vectors
$$ w_1=(-b/c,a/c,0) \ \text{and} \ w_2=(0,-c/a,b/a).$$

Let $s \in {\Bbb F}_q$ be such that
$$ e_1=w_1+sw_2$$ satisfies $ ||e_1||=1.$ Such $s$ exists, by the Lagrange theorem on quadratic forms (or can be verified by direct calculation).

Consider now a six-vector $v_3=[e_1,w_1]$. By construction, $v_3$ is orthogonal to both
$v_1$ and $v_2$. It is also a null vector, as $e_1\cdot e_1=1$, while $w_1\cdot w_1=-1$.

This completes the proof of Lemma \ref{nullvectors}.

\vskip.125in

Let $d=2m+1$, then from the above lemma there are $m$ mutually orthogonal
null
vectors $v_1,\dots,v_m$, such that their $d$th coordinate is zero.
Now let $A\subset \mathbb F_q$ be an arithmetic progression of length
$n$ and $u=(0,\dots,0,1)$.
Consider the set
$$E=\{t_i v_i + a u  \ \text{for}\  i=1,\dots,m: t_i\in \mathbb F_q, a
\in A\}.$$

We have
$$ |E|=q^m \cdot |A|=q^m \cdot n.$$

For any $x,y\in E$ we have from orthogonality that
$$\|x-y\| = \|t_1u_1+av-t_2u_2-a'v_2\|=(a-a')^2,$$
so $|\Delta(E)|\leq 2n-1$.

It follows that if we choose $2n=cq$, we have constructed, for any
small $c$, a set $E$ of
$\frac{1}{2} cq^{\frac{d+1}{2}}$ generating fewer than $cq$ distances.
This completes the proof in the case $d \ge 5$.

If $d=3$, and $-1$ is a square, take the null vector $v=(1,i,0)$ and
$u=(0,0,1)$. If $-1$ is not a square, take the null vector $v=(a,b,c)$
such that no entry can be zero, and let $u=(-b,a,0)$. The proof then
proceeds as above.

\vskip.125in

\subsection{Proof of Theorem \ref{sphericalerdos}, claim (i)}
Since $E \subset S$,
$$\|x-y\|=(x-y) \cdot (x-y)=2-2x \cdot y,$$ so counting distances on
the sphere is the same as counting dot products.

Since now $E$ is a subset of the
sphere, it does not contain the origin and
$$ |E \cap l_k| \leq 2.$$

Thus we conclude from the estimate (\ref{L2}) of Theorem \ref{geom} that
\begin{equation} \label{spherepizdets} \begin{array}{lll}|E|q^{2d-1}
\sum_{k \neq(0,\ldots,0)} |E \cap l_k| {|\widehat{E}(k)|}^2 &\leq&
2|E|q^{2d-1} \sum_{k \neq(0,\ldots,0)} {|\widehat{E}(k)|}^2\\ \hfill
\\
& \leq & 2|E|q^{2d-1} q^{-d} \sum_x
E^2(x)\\ \hfill \\
&=&2{|E|}^2q^{d-1}.\end{array}
\end{equation}
By application of the Cauchy-Schwartz inequality (\ref{csagain})
we conclude that if $|E| \ge cq^{\frac{d}{2}},$ we have
$$ |\Delta(E)| \ge Cq.$$
\subsection{Proof of the claim (iv)}  The claim for $x\cdot y \neq 0$
follows immediately from Corollary \ref{L2dot}, without requiring $d$
to be odd. The case $x\cdot y=0$ will be addressed further in Section
\ref{long}.

\vskip.125in

\vskip.125in

\subsection{Proof of Theorem \ref{sphericalerdos}, claim (ii)} \label{long}
We now turn to the proof of (\ref{spherepointwise}). We will not
distinguish between even and odd $d$ until it becomes necessary. We
proceed as in (\ref{break}) in the proof of Corollary
(\ref{pointwise}) by writing
$$ \nu(t)={|E|}^2q^{-1}+R(t),$$ and apply the Cauchy-Schwartz
inequality to $R^2(t)$. This time, however, instead of dominating the
sum over $E$ by the sum over ${\mathbb F}_q^d$, we dominate the sum
over $E$ by the sum over the sphere $S$ using the assumption that $E
\subset S$. This yields
\begin{equation}\begin{array}{lll}R^2(t) &\leq& q^{-2}|E| \sum_{x \in
S}  \sum_{s,s' \not=0}
\sum_{y,y' \in E}
\chi(sx \cdot y-s'x \cdot y') \chi(t(s'-s)) \\ \hfill \\
&=& q^{d-2}|E| \sum_{s,s' \not=0} \chi(t(s'-s))
\sum_{y,y' \in E} \widehat{S}(s'y'-sy) \\ \hfill \\ &=&I+II,\end{array}
\label{split}\end{equation}
where the term $I$ corresponds to the case $s'y'=sy$. One of the keys
to this argument is that since $E$ is a subset of the sphere,
$sy=s'y'$ can only happen if $y= \pm y'$ and $s= \pm s'$.

Lemma (\ref{sphereft}) below tells us that $\widehat{S}(0, \dots,
0)=q^{-1}+$ lower order terms (unless $d=2$), and it follows that
\begin{equation}
I\leq q^{d-2}|E|^2.
\label{Iok}\end{equation}
To estimate the term $II$, we have to use the explicit form of the
Fourier transform of the discrete sphere. For the reader's convenience
we replicate one of the arguments in \cite{IR07}.

\begin{lemma} \label{sphereft}Let
$$
S_r=\{(x_1,\ldots,x_d)\in{\mathbb F}_q^d:\,x_1^2+\ldots+x_d^2=r\},
$$
Then for $k\in\mathbb F_q^d$,
\begin{equation}
\widehat S_r(k)= q^{-1}\delta(k)+K^d
q^{-\frac{d+2}{2}}\sum_{j\in{\mathbb
F}_q^*}\chi\left(\frac{\|k\|}{4j}+rj\right)\eta^d(-j),
\label{fts}
\end{equation}
where the notation $\delta(k)=1$ if $k=(0\ldots,0)$ and $\delta(k)=0$
otherwise. The constant $K$ equals $\pm1$ or $\pm i$, depending on
$q$, and $\eta$ is the quadratic multiplicative character (or the
Legendre symbol) of ${\mathbb F}_q^*$.
\end{lemma}

\subsubsection{Proof of Lemma \ref{sphereft}} For any $k\in{\mathbb
F}^d_q$, we have
\begin{equation} \label{sphereparade}
\begin{array}{llllll} \widehat{S}_r(k)&=&
q^{-d} \sum_{x \in {\mathbb F}^d_q} q^{-1} \sum_{j \in {\mathbb F}_q} \chi(
j(\|x\|-r)) \chi( -  x \cdot k)\\ \hfill \\&=&q^{-1}\delta(k) +
q^{-d-1} \sum_{j \in {\mathbb F}^{*}_q} \chi(-jr) \sum_{x}
\chi( j\|x\|) \chi(- x \cdot k)\\ \hfill
\\&=&q^{-1}\delta(k)+   K^d q^{-\frac{d+2}{2}} \sum_{j \in {\mathbb
F}^{*}_q}
\chi\left(\frac{\|k\|}{4j}+jr\right)\eta^d(-j).\end{array}\end{equation}
In the line before last we have completed the square, changed $j$ to
$-j$, and used $d$ times the Gauss sum
\begin{equation}
\sum_{c\in {\mathbb F}_q} \chi(jc^2) = \eta(j)\sum_{c\in{\mathbb
F}_q}\eta(c)\chi(c)=\eta(j)\sum_{c\in{\mathbb F}_q^*}\eta(c)\chi(c) =
K\sqrt{q}\,\eta(j),
\label{gauss}\end{equation} where $K=\pm i$ or $\pm 1$, depending on
$q$ and $\eta(0)=0$. See any standard text on finite fields for
background and basic results about Gauss sums. Note that we have
assumed that $\chi=\chi_1$ is the principal additive character of the
field ${\mathbb F}_q$ (which means that for $t\in{\mathbb F}_q$, and
$q=p^s$, where $p$ is a prime, $\chi(t)=e^{\frac{2\pi i {\rm
Tr}(t)}{p}}$, where ${\rm Tr}:\,{\mathbb F}_q\mapsto{\mathbb F}_p$ is
the principal trace, see e.g. \cite{NL}.) The specific choice of a
principal character is of no consequence to the calculations in this
paper.

We remark that for even $d$, the sum in the last line of
(\ref{sphereparade}) is the Kloosterman sum, while for odd $d$ the
presence of the quadratic character $\eta$ would reduce it via the
Gauss sum to a ``cosine'', which is nonzero only if
$\theta^2\equiv\frac{ r\|k\|}{4}$ is a square in ${\mathbb F}_q^*$, in
which case
\begin{equation}
\sum_{j \in {\mathbb F}^{*}_q}
\chi\left(\frac{\|k\|}{4j}+jr\right)\eta(-j) =
K\sqrt{q}\,\eta(-\|k\|^2) (\chi(2\theta)+\chi(-2\theta)).
\label{cosine}\end{equation}

\medskip
We now return to the proof of (\ref{spherepointwise}). From now on,
let $K,K'$ stand for complex numbers of modulus $1$ that may change
from line to line. We now continue with the estimation of the term
$II$ in (\ref{split}). Namely, we have
$$ II= q^{d-2}|E|\sum_{y,y' \in E} \sum_{s,s'\in{\mathbb F}_q^*,
s'y'\neq sy}\widehat{S}(s'y'-sy) \chi(t(s'-s))= III + IIII,$$
where the term $III$ corresponds to the case the case $y=y'$, when $s\neq
s'$.
Then we have
$$III= q^{d-2}|E| \sum_{y \in E}\sum_{s,s' \in{\mathbb F}_q^*, s\neq
s'} \widehat{S}((s'-s)y) \chi(t(s'-s)).$$

Observe that $s'-s$ runs through each value in ${\mathbb F}_q^*$
exactly $q-1$ times. Also, $||y||=1$ since $E\subset S$. Therefore,
using Lemma \ref{fts}, we have
\begin{equation}\label{three}\begin{array}{lll}
III &=& K q^{d-2}|E| \sum_{y \in E} (q-1)
q^{-\frac{d+2}{2}}\sum_{s,j\in{\mathbb F}_q^*}
\chi\left(\frac{s^2}{4j}+ ts +j\right)\eta^d(j)
\\ \hfill \\
&=& K q^{d-2}|E| \sum_{y \in E} (q-1)\cdot
q^{-\frac{d+2}{2}}\sum_{s,j\in{\mathbb F}_q^*} \chi\left(\frac{(s+
2jt)^2}{4j}- jt^2 +j\right)\eta^d(j)
\\ \hfill \\
&=&K \frac{q-1}{q} q^{\frac{d-4}{2}}|E| \sum_{y \in E}
\sum_{j\in{\mathbb F}_q^*} \chi(j-jt^2)\eta^d(j)
[-\chi(jt^2)+K'\sqrt{q}\eta(j)],
\end{array}
\end{equation}
where the last line follows by (\ref{gauss}).

We now consider the case $t^2=\|y\|=1$. We have
$$ III_{t^2=1} \approx q^{\frac{d-4}{2}}|E| \sum_{y \in E}
\sum_{j\in{\mathbb F}_q^*} \eta^d(j) [-\chi(j)+K \sqrt{q}\eta(j)]. $$

Since $\sum_{j\in{\mathbb F}_q^*}\eta(j)=0$, the worst case scenario
is when $d$ is odd. Then the summation in $j$ contributes an extra
factor $q-1$ to $K\sqrt{q}$ in the last bracket. If $d$ is even then
the summation in $j$ is the Gauss sum, which is smaller by the factor
of $\sqrt{q}$. In either case, we have
\begin{equation} |III_{t^2=1}| \leq 2q^{\frac{d-1}{2}}|E|^2.
\label{hmodd}\end{equation}

If $t^2\neq 1$, the estimate (\ref{hmodd}) improves by factor
$\sqrt{q}$, as the worst case scenario is now when $d$ is even, and it
only contributes a Gauss sum to the term $K\sqrt{q}$:
\begin{equation}\label{analysis}\sum_{j\in{\mathbb F}_q^*}
\chi(j-jt^2)\eta^d(j) [-\chi(jt^2)+K\sqrt{q}\eta(j)].\end{equation}
Observe however, that in either case, for $d\geq2$ the estimate for
the term $III$ is majorated by (\ref{Iok}).

\vskip.125in

Finally, let us consider the term $IIII$:
\begin{equation} IIII = q^{d-2}|E| \sum_{y,y' \in E, y\neq
y'}\sum_{s,s' \in{\mathbb F}_q^*}
\widehat{S}(s'y'-sy)\chi[t(s'-s)]\label{four}\end{equation}
Our goal is to prove the following estimate:

\begin{equation}\label{prove}
|IIII=IIII(t)|\;\;\lesssim \;\;q^{\frac{d-4}{2}}|E|^3 +
q^{\frac{d-2}{2}}|E| \sup_{\tau\in{\mathbb F}_q}|R(\tau)|.
\end{equation}
and we are able to do it only for even values of $d$. (For odd $d$ the
estimate will be definitely worse by $\sqrt{q}$ for $t^2=1$ and seems
to be highly non-trivial for other values of $t$, see (\ref{oscil})
below.) Note that we can always write
$\sup_{\tau\in{\mathbb F}_q} \nu(\tau)$ instead of
$\sup_{\tau\in{\mathbb F}_q}|R(\tau)|,$ as the regular term
$\frac{|E|^2}{q}$ can be absorbed into the first term in
(\ref{prove}).

\vskip.125in

We verify (\ref{prove}) below and will now show how it suffices to
complete the proof of (\ref{spherepointwise}). Indeed, assuming
(\ref{prove}) and bringing in the estimate (\ref{Iok}), which
dominated the terms $I,\,III$, we conclude that for all $t$,
$$ R^2(t) \;\lesssim \;q^{d-2} |E|^2 + q^{\frac{d-4}{2}}|E|^3 +
q^{\frac{d-2}{2}}|E| \sup_{\tau\in{\mathbb F}_q} |R(\tau)|,$$ which
implies that  the same estimate holds for $\sup_{\tau\in{\mathbb
F}_q}R^2(\tau)$.

Assuming that for some large enough $C$, we have $Cq^{\frac{d}{2}}
\leq |E|$ clearly implies that now
$$ |R(t)|\;\leq\;\frac{100}{\sqrt{C}}\frac{|E|^2}{q},\;\;\;\forall
t\in{\mathbb F}_q, $$
where the constant $100$ is basically to majorate the number of cases
that has been considered. For odd $d$ the last two terms in the latter
estimate for $R$ are worse by the factor $\sqrt{q}$ which implies the
estimate (\ref{remainder}) for all $t$, thus the claim (ii) of Theorem
\ref{sphericalerdos}. As for even $d$, every dot product $t\in
{\mathbb F}_q$ occurs and the claim (iv) of Theorem
\ref{sphericalerdos} follows, provided that we can demonstrate
(\ref{prove}).

\subsubsection{Finale of the proof of claim (ii) -- the estimate
(\ref{prove})} \label{finale}
In the estimates that follow we write
$$ \sum_{y,y'}\;\;\;\;\mbox{ instead of }\;\;\;\;\sum_{y,y' \in E, y\neq
y'}.$$

Let us first extend the summation in (\ref{four}) from $s'\in{\mathbb
F}_q^*$ to $s'\in{\mathbb F}_q$. If we do so, it follows from Lemma
\ref{fts} that we pick up the following term $T$ to $IIII$:
$$\begin{array}{lll}
T &=& q^{d-2}|E| \sum_{y,y' }\sum_{s\in{\mathbb F}_q^*}
\widehat{S}(sy)\chi(ts) \\ \hfill \\
&=& K q^{\frac{d-6}{2}}|E| \sum_{y,y'} \sum_{s,j\in{\mathbb F}_q^*}
\chi\left(\frac{s^2}{4j}+ts+j\right)\eta^d(j)
 \\ \hfill \\ &=&
K q^{\frac{d-6}{2}}|E| \sum_{y,y'} \sum_{s,j\in{\mathbb F}_q^*}
\chi\left(\frac{(s+ 2jt)^2}{4j}- jt^2 +j\right) \eta^d(j)
 \\ \hfill \\ &=&
 K q^{\frac{d-6}{2}}|E| \sum_{y,y'} \sum_{j\in{\mathbb
F}_q^*}\chi(j-jt^2)\eta^d(j) [-\chi(jt^2)+K'\sqrt{q}\eta(j)]
 \end{array}$$
 The analysis of the summation in $j$ now in essence replicates that
for the term $III$, see (\ref{three}--\ref{analysis}). If $t^2\neq 1$
and  $d$ is even, using the Gauss sum formula (\ref{gauss}) we obtain
 \begin{equation} \label{4ok}
|T| \;\leq \;2q^{\frac{d-4}{2}}|E|^3,\end{equation}
 which improves by factor $\sqrt{q}$ if $d$ is odd. If $t^2=1$, for
even $d$, the term $T$ satisfies a better (by a factor $q$) estimate
than (\ref{4ok}). However, for odd $d$ we would only get
 $|T| \leq q^{\frac{d-3}{2}}|E|^3$, which would not give an
improvement over the bounds we already have we already have
in (\ref{remainder}). Hence, up to now, the only case we are not able
to handle is odd $d$ and $t^2=1$.

\vskip.125in

Thus we will further attempt to establish (\ref{prove}) for the
quantity $X$, which equals $IIII$, wherein the summation in $s'$ has
been extended over the whole field ${\mathbb F}_q$. Using Lemma
\ref{fts} we have, after changing $s'$ to $-s'$,  and using
$\|y\|=\|y'\|=1$:
\begin{equation}\label{ox} X =K q^{\frac{d-6}{2}}|E|
\sum_{y,y'} \sum_{s,j\in{\mathbb F}_q^*,s'\in{\mathbb F}_q}\eta^d(j)
\chi\left(\frac{s^2+2(y\cdot y')ss' + {s'}^2+4tj(s+s')} {4j}+
j\right) \end{equation}

We complete the square under $\chi$ as follows
$$ s^2+2(y\cdot y')s s' + {s'}^2+4tj(s+s') = (s+s')^2 + 4tj(s+s') +
2\alpha s'(s+s'-s')$$ where $\alpha=\alpha(y,y') = y\cdot y'-1$, and
we shall further analyze the possibilities $\alpha\neq 0,-2$
separately: they occur when $y\cdot y'=\pm1$, respectively.

We rewrite the latter quadratic form as
$$ [(s+s') + (2tj+\alpha s')]^2 - 2\alpha s'^2 - (2tj+\alpha s')^2. $$

We now have a new variable $c = (s+s') + (2tj+\alpha s')$, which is in
${\mathbb F}_q$. Since (\ref{ox}) is symmetric with respect to $s$ and
$s'$, we can assume that, in fact, $s\in {\mathbb F}_q,\,s'\in{\mathbb
F}_q^*$, so for each $s',j$ the change $s\mapsto c$ is non-degenerate.
Changing the notation from $-s'$ to $s$ we therefore have, using the
Gauss sum formula
\begin{equation}\label{ex}\begin{array}{lll}
 X
&=& K q^{\frac{d-6}{2}}|E| \sum_{y,y'} \sum_{s,j\in{\mathbb
F}_q^*,c\in{\mathbb F}_q}\eta^d(j)
\chi\left(\frac{c^2 - (2\alpha+\alpha^2)s^2}{4j} + t\alpha s + j(1 -
t^2) \right)
 \\ \hfill \\ &=&
 X_1+X_{-1}+X',
 \end{array}\end{equation}
 where $X_1$ has only summation  in $y,y'$ such that $y\cdot y'=1$
($\alpha=0$), $X_{-1}$ has only summation  in $y,y'$ such that $y\cdot
y'=-1$ ($\alpha=-2$), and $X'$ includes the rest of $y,y'\in E$.

Observe that in either case we already have a Gauss sum in $c$, so we write
\begin{equation}\label{exp}
X'=K q^{\frac{d-5}{2}}|E| \sum_{y\cdot y'\neq\pm1}
\sum_{s,j\in{\mathbb F}_q^*}\eta^{d+1}(j)\chi\left( \frac{- a \left(s
-\frac{ 2jt\alpha}{a}\right)^2}{4j} +
j\left(\frac{t^2\alpha}{2+\alpha}+(1-t^2)\right) \right),
\end{equation} provided that $a=2\alpha+\alpha^2\neq0$.

Before we proceed with the main term $X'$, let us deal with the cases
$\alpha=0,-2$ which would make the completion of the square in the
transition from (\ref{ex}) to (\ref{exp}) incorrect.

If $\alpha=0$, we confront the sum
$$ X_1= K q^{\frac{d-5}{2}}|E| \sum_{y\cdot y'=1} \sum_{s,j\in{\mathbb
F}_q^*}\eta^{d+1}(j)
\chi( j(1 - t^2)). $$

If $d$ is even, the worst case scenario is $t^2\neq1$, when the sum in
$s$ and Gauss sum in $j$ contribute the factor $q^{3/2}$. Hence
\begin{equation} |X_{1}|\,\leq \,2q^{\frac{d-2}{2}}|E| \sup_\tau
\nu(\tau),\;\;\mbox{ for even } d, \label{good}\end{equation}
in accordance with (\ref{prove}). If $d$ is odd, the same, or in fact,
better bound holds unless $t^2=1$, when (\ref{good}) gets worse by
factor $\sqrt{q}$.

If $\alpha=-2$, we analyze the sum
$$ X_{-1}= K q^{\frac{d-5}{2}}|E| \sum_{y\cdot y'=-1}
\sum_{s,j\in{\mathbb F}_q^*}\eta^{d+1}(j)
\chi( j(1 - t^2) - 2ts). $$
If $d$ is even, $X_{-1}$ is still bounded by (\ref{good}) -- the worst
case scenario now is $t=0$; if $d$ is odd, the bound is better than
(\ref{good}) by factor $\sqrt{q}$.

\vskip.125in

Finally, we turn to $X'$, the case $a\neq0$, and once again, the only
situation we have not been able to handle so far is $d$ odd and
$t^2=1$.

Now taking advantage of the Gauss sum in $s$ in (\ref{exp}) we have
$$\begin{array}{lll}
 X'&=&
K q^{\frac{d-5}{2}}|E| \sum_{y\cdot y'\neq\pm1} \sum_{j\in{\mathbb
F}_q^*}\eta^{d+1}(j)\chi\left(j\left(\frac{t^2\alpha}{2+\alpha}+(1-t^2)\right)\right)
\\ \hfill \\
&&\hspace{1in}\times\;\;\left[-\chi\left(
-j\frac{t^2\alpha}{2+\alpha}\right) +K'\sqrt{q}\eta(a)\eta(j)\right]\\
\hfil \\
&=&X'_1+X'_2,\end{array}$$
according to the two terms in the last bracket.

We have
$$ X'_1 = K q^{\frac{d-5}{2}}|E| \sum_{y\cdot y'\neq\pm1}
\sum_{j\in{\mathbb F}_q^*}\eta^{d+1}(j)\chi(j(1-t^2)). $$ For even
$d$, the worst case scenario occurs when $t^2\neq 1$, the Gauss sum in
$j$ then leads to $X_1'$ to be dominated by the first term in
(\ref{prove}). The latter bound will get worse by factor $\sqrt{q}$
only if $d$ is odd and $t^2=1$. For the quantity $X_2'$ we obtain:
\begin{equation}\begin{array}{lll}X_2'&=&
K q^{\frac{d-4}{2}}|E| \sum_{y\cdot y'\neq\pm1} \eta(a)
\sum_{j\in{\mathbb
F}_q^*}\eta^{d}(j)\chi\left(\left(1-\frac{2}{2+\alpha}t^2\right)j\right)\\
\hfill \\
&=&
K q^{\frac{d-4}{2}}|E| \sum_{y\cdot y'\neq\pm1} \eta[(y\cdot y')^2-1]
\sum_{j\in{\mathbb F}_q^*}\eta^{d}(j)\chi\left(\left(1-\frac{2}{y\cdot
y'+1}t^2\right)j\right).\end{array}\label{X2}\end{equation}

There are two cases here: $y\cdot y'=2t^2-1$ and otherwise. First
consider the latter case. Then if $d$ is even, $X'_2$, subject to this
extra constraint, satisfies the estimate (\ref{prove}), as the
summation in $j$ simply yields $-1$.
If $d$ is odd, however, there is a major problem, as then we have
\begin{equation}\label{oscil}\begin{array}{lrr}
& q^{\frac{d-4}{2}}|E| \sum_{y\cdot y'\neq 2t^2-1,\pm1} \eta((y\cdot
y')^2-1) \sum_{j\in{\mathbb
F}_q^*}\eta^{d}(j)\chi\left(\left(1-\frac{2}{y\cdot
y'+1}t^2\right)j\right) \\ \hfill \\
&=Kq^{\frac{d-3}{2}}|E| \sum_{y\cdot y'\neq 2t^2-1,\pm1} \eta((y\cdot
y')-1) \eta((y\cdot y')+1-2t^2), \end{array}. \end{equation}

It follows that to improve on the trivial bound
$q^{\frac{d-3}{2}}|E|^3$ one would have to establish a cancelation in
the  multiplicative character sum in (\ref{oscil}).

We finish by adding the constraint $y\cdot y'=2t^2-1$ to $X_2'$ in
(\ref{X2}). Dealing with this does not represent any difficulty. For
even $d$ we have
$$ q^{\frac{d-4}{2}}|E| \left|\sum_{\pm1\neq y\cdot y'=2t^2-1}
\eta[(y\cdot y')^2-1] \sum_{j\in{\mathbb F}_q^*}\eta^{d}(j) \right|
\leq q^{\frac{d-2}{2}}|E|\sup_{\tau\in{\mathbb F}_q} \nu(\tau), $$
 and zero in the right-hand side for odd $d$. This proves
(\ref{prove}) and (\ref{spherepointwise}) follows.

\vskip.125in

\subsection{Proof of Theorem \ref{sphericalerdos}, optimality claims
(iii) and (v)}
\label{nodot}
\vskip.125in

We establish (\ref{oddcounter}) as the estimate (\ref{evencounter})
follows immediately from the same construction.

\subsection{Construction in the case $d \not=5$:} Suppose that
${\mathbb F}_q$ does not contain
$i=\sqrt{-1}$. Let
$$ S^2=\{x \in {\mathbb F}_q^3: x_1^2+x_2^2+x_3^2=1\},$$ and let
${Z_2}$ denote the maximal subset of $S^2$ such that ${Z_2} \cap
(-{Z_2})=\emptyset$. Then if $u,v \in S^2$, then $u \cdot v=-1$ if and
only if $u=-v$. To see this, without loss of generality let
$v=(0,0,1)$. Then the condition
$$ u \cdot v=-1$$ reduces to
$$ u_3=-1,$$ and
\begin{equation} \label{zhopa} u_1^2+u_2^2=0. \end{equation}

Since, by assumption, ${\mathbb F}_q$ does not contain $\sqrt{-1}$,
(\ref{zhopa}) can only happen if $u_1=u_2=0$, and so $u=-v$. Since
${Z_2} \cap (-{Z_2})=\emptyset$, the condition $u \cdot v=-1$ in
${Z_2}$ is never satisfied.

Let $d=2k+1$ with $k \ge 3$. Let $H$ denote sub-space of ${\Bbb
F}_q^{2k-2}$ generated by the mutually orthogonal null-vectors given
by Lemma \ref{nullvectors}. Let
$$ E=Z_2 \times H.$$

It follows that
$$ |E| \approx q^2 \cdot q^{k-1}=q^{k+1}=q^{\frac{d+1}{2}}.$$

Let $(x',x'')$ and $(y',y'')$ be elements of $E$. Then
$$ (x',x'') \cdot (y', y'')=x' \cdot y' \not=-1.$$

Moreover,
$$ ||(x',x'')||=||x'||+||x''||=||x'||=1,$$ so $E \subset S^{2k}$ where
$$ S^{2k}=\{x \in {\Bbb F}_q^{2k+1}: x_1^2+\dots+x_{2k+1}^2=1\}.$$

This completes the construction in the case $d \not=5$.

\subsection{Construction in the case $d=5$} Let
$$ u=(a,b,c,0,0) \ \text{where} \ a^2+b^2+c^2=0.$$

Let
$$ v=(-b/c,a/c,0,0,0) \ \text{and} \ w=(0,-c/a,b/a,0,0).$$

Let $s \in {\Bbb F}_q$ be such that
$$ e=v+sw$$ satisfies
$$ ||e||=c^2 \ \text{for some} \ c \in {\Bbb F}_q^{*}.$$

The existence of such a $c$ is verified by a direct calculation. Now
let $e'=\frac{e}{c}$, which results in $||e'||=1$.

Observe by a direct calculation that
$$ u \cdot e=0 \ \text{for all} \ s \in {\Bbb F}_q.$$

Let $Z_2$ be as above and let $O$ denote the orthogonal transformation
that maps
$$ \{(x_1,x_2,x_3,0,0): x_j \in {\Bbb F}_q\}$$ to the three
dimensional sub-space of ${\Bbb F}_q^5$ spanned by $e'$, $(0,0,0,1,0)$
and $(0,0,0,0,1)$. Let $Z'_2$ denote the image of $Z_2$ under $O$.

Define
$$ E=\{tu+Z'_2: t \in {\Bbb F}_q\}.$$

Then $|E| \approx q^3$ and for any $t,t' \in {\Bbb F}_q$ and $z,z' \in Z'_2$,
$$ (tu+z) \cdot (t'u+z')=tt'u \cdot u+tu \cdot z'+t'u \cdot z+z \cdot z'$$
$$=z \cdot z' \not=-1$$ by construction. This completes the
construction in the case $d=5$.

\vskip.125in

\subsection{Proof of the conditionally optimal result (Theorem
\ref{uniformsphere})}

\vskip.125in

Once again we use the estimate (\ref{L2}) which tells us that
$$ \sum_t \nu^2(t) \leq {|E|}^4q^{-1}+|E|q^{2d-1} \sum_{k
\neq(0,\ldots,0)} |E \cap l_k|
{|\widehat{E}(k)|}^2.$$

Now,
$$ |E|q^{2d-1} \sum_{k \neq(0,\ldots,0)} |E \cap l_k|
{|\widehat{E}(k)|}^2 \leq 2|E|q^{2d-1} \sum_{k\in C(E)}
{|\widehat{E}(k)|}^2,$$ where
$$C(E)=\bigcup_{t \in {\mathbb F}_q} tE.$$
Furthermore,
\begin{equation} \label{plugunif} \begin{array}{lll} |E|q^{2d-1}
\sum_{k\in C(E)} {|\widehat{E}(k)|}^2&=& |E|q^{2d-1} q^{-2d}
\sum_{y,y' \in E} \sum_{k\in C(E)} \chi((y'-y) \cdot k)\\ \hfill \\
&=&|E|q^{-1} \sum_{y,y' \in E} \sum_t \sum_{x \in E} \chi((y-y') \cdot tx)$$
\\ \hfill \\
&=&{|E|}^3+|E| \sum_{(y-y') \cdot x=0; y \not=y'}
E(x)E(y)E(y'). \end{array}\end{equation}

Since $E$ is assumed to be uniformly distributed,
$$ \sum_{(y-y') \cdot x=0} E(x) \leq C|E|q^{-1},$$ plugging this into
(\ref{plugunif}) we obtain
${|E|}^4q^{-1}$. Using (\ref{csagain}) once again we complete the
proof. Observe that the assumption that $|E| \ge Cq$ is implicit in
the uniform distributivity assumption.


\vskip.125in

\newpage

\begin{thebibliography}{17}

\vskip.125in

\bibitem{AK97} Alon and Krivelevich, {\it Constructive bounds for a
Ramsey-type problem}, Graphs and Combinatorics \textbf{13} (1997),
217-225.

\bibitem{BGK06} J. Bourgain, A. A. Glibichuk and S. V. Konyagin. {\it
Estimates for the number
of sums and products and for exponential sums in fields of prime order}. J.
London Math. Soc. (2) \textbf{73} (2006), 380--398.

\bibitem{BKT04} J. Bourgain, N. Katz and T. Tao. {\it A sum-product estimate
in finite fields, and applications}. Geom. Func. Anal. \textbf{14} (2004),
27--57.

\bibitem{C04} E. Croot. {\it Sums of the Form $1/x_1^k+\dots 1/x_n^k$ modulo
a prime}. Integers \textbf{4} (2004).

\bibitem{Erd05} B. Erdo\~{g}an. {\it A bilinear Fourier extension
theorem and applications to the
distance set problem.}  Int. Math. Res. Not. {\bf 23} (2005), 1411--1425.

\bibitem{Erd45} P. Erd\H os. {\it On sets of distances of n points.}
Amer. Math. Monthly. \textbf{53} (1946), 248--250.


\bibitem{G07} M. Garaev. {\it The sum-product estimate for large subsets of
prime fields}. Preprint, 2007.

\bibitem{G06} A. A. Glibichuk. {\it Combinatorial properties of sets of
residues modulo a
prime and the Erd\H os-Graham problem}. Mat. Zametki \textbf{79} (2006),
384--395;
translation in: Math. Notes \textbf{79} (2006), 356--365.

\bibitem{GK06} A. Glibichuk and S. Konyagin. {\it Additive properties of
product sets in fields of prime order}. Centre de Recherches Mathematiques,
Proceedings and Lecture Notes, 2006.

\bibitem{HI07} D. Hart and A. Iosevich. {\it Sums and products in
finite fields: an integral geometric viewpoint}. Preprint, arxiv.org,
2007.

\bibitem{HIS07} D. Hart, A. Iosevich and J. Solymosi. {\it Sums and products
in finite fields via Kloosterman sums}. Int. Math. Res. Not.  To appear,
2007.

\bibitem{IR07} A. Iosevich and M. Rudnev. {\it Erd\H os distance problem in
vector spaces over finite fields}. Trans. Amer. Math. Soc. To appear, 2007.

\bibitem{IRU07} A. Iosevich, M. Rudnev and I. Uriarte-Tuero. {\it
Theory of dimension for large discrete sets and applications.}
Preprint, arxiv.org, 2007.

\bibitem{KS07}  Nets Hawk Katz and Chun-Yen Shen. {\it Garaev's Inequality
in finite fields not of prime order}. Preprint, arxiv.org,
2007.

\bibitem{NL} R. Lidl and H. Niederrieter. {\em Finite Fields.}
Encyclopedia of Mathematics and its Applications {\bf 20},
Addison--Wesley 1983.

\bibitem{Ma02} J. Matousek. {\it Lectures on Discrete Geometry,}
Graduate Texts in Mathematics. Springer \textbf{202}, 2002.

\bibitem{Mat87}
P. Mattila.
{\it Spherical averages of {F}ourier transforms of measures with finite
 energy; dimension of intersections and distance sets.}
Mathematika {\bf 34}(2) (1987), 207--228.

\bibitem{St93} E. Stein. {\it Harmonic Analysis}. Princeton University
Press, 1993.

\bibitem{TV06} T. Tao and V. Vu. {\it Additive Combinatorics}. Cambridge
University Press, 2006.

\bibitem{V07} V. Vu. {\it Sum-Product estimates via directed
expanders}. Preprint, arxiv.org,
2007.

\bibitem{We48} A. Weil. {\it On some exponential sums}. Proc. Nat.
Acad. Sci. U.S.A. \textbf{34} (1948), 204--207.

\end{thebibliography}
\end{document}